\documentclass[12pt,a4paper]{article}
\usepackage[utf8]{inputenc}
\newcommand{\pfrac}[2]{\frac{\partial{#1}}{\partial{#2}}}
\usepackage[margin=2cm]{geometry}
\usepackage{bm} % for bold maths symbols
\usepackage{enumitem} % for more flexibility with lists
\usepackage{array}  % to give vertical centering in tables
\usepackage[pdftex]{graphicx} % for including images
\usepackage{subcaption} % for side by side figures
\usepackage{caption}
\usepackage{amssymb} % for mathbb
\usepackage{amsmath} % for align and tfrac
\usepackage{natbib}
\usepackage{hyperref} % for interactable links
\usepackage{tikz-cd}
\usepackage{lineno}
\usepackage{float}
\usepackage{multicol}
\usepackage{algorithm}
\usepackage{algpseudocode}
\usepackage{authblk} % for author and affiliations
\usepackage{fancyhdr} % For headers and footers
\newcommand{\grad}[1]{{\bm{\nabla}{#1}}}
\newcommand{\dv}[1]{\bm{\nabla\cdot{#1}}}

\DeclareMathOperator{\diag}{diag}

\begin{document}

\pagestyle{fancy}
\fancyhf{}
\lfoot{\textcopyright \ Crown Copyright, Met Office}
\cfoot{\thepage}

\title{Fast-wave slow-wave spectral deferred correction methods applied to the compressible Euler equations}

\author[1,2]{Alex Brown}
\author[3]{Joscha Fregin}
\author[1]{Thomas Bendall}
\author[1]{Thomas Melvin}
\author[3]{Daniel Ruprecht}
\author[2]{Jemma Shipton}

\affil[1]{Dynamics Research, Met Office, Exeter, UK}
\affil[2]{Department of Mathematics, University of Exeter, UK}
\affil[3]{Chair Computational Mathematics, Institute of Mathematics, Hamburg University of Technology, Hamburg, Germany}
\date{}

\maketitle

\begin{abstract}
This paper investigates the application of a fast-wave slow-wave spectral deferred correction time-stepping method (FWSW-SDC) to the compressible Euler equations.
The resulting model achieves arbitrary order accuracy in time, demonstrating robust performance in standard benchmark idealised test cases for dynamical cores used for numerical weather prediction. 
The model uses a compatible finite element spatial discretisation, achieving good linear wave dispersion properties without spurious computational modes. 
A convergence test confirms the model's high temporal accuracy. Arbitrarily high spatial-temporal convergence is demonstrated using a gravity wave test case. 
The model is further extended to include the parametrisation of a simple physics process by adding two phases of moisture and its validity is demonstrated for a rising thermal problem.
Finally, a baroclinic wave in simulated in a Cartesian domain.
\end{abstract}

\section{Introduction} \label{sec:intro}

\thispagestyle{fancy}  
Semi-implicit timestepping methods are often used in numerical weather prediction (NWP) and climate modelling, see~\cite{robert1981stable}. 
They allow for large time steps to be taken, stabilising the simulation of fast-moving waves, whilst maintaining the accurate representation of the slow-moving waves. 
They are computationally more efficient than a fully-implicit method, and can maintain accuracy when taking large time steps by not damping the slow-moving waves. 
Semi-implicit methods, such as the one used in the Met Office's next dynamical core GungHo, see~\cite{melvin2024mixed,melvin2019mixed}, often use an off-centered time integration scheme to treat the implicit terms. Such schemes are generally only second order when fully centered, and can reduce to first order when off-centered. Other time-split methods commonly used in atmospheric modeling are split-explicit methods, such as those used in MPAS (see~\cite{skamarock2012multiscale}) and WRF (see~\cite{skamarock2008description}). They use different explicit methods and timestep sizes for terms responsible for fast and slow moving waves and use a smaller timestep than semi-implicit methods, but sometimes use higher order timestepping methods.

Spectral deferred correction methods (SDC), originally introduced by~\cite{dutt2000spectral}, have been used in a wide range of fields. 
Examples are climate modelling by~\cite{jia2013spectral}, particle modelling by~\cite{winkel2015high} and phase-field problems by~\cite{feng2015long}. 
\cite{OngSpiteri2020}~provide a review of variants and applications of deferred correction methods.
SDC are arbitrarily high order time integrators, which are competitive in terms of accuracy, stability and cost to existing time integrators such as Runge-Kutta methods, see~\cite{ruprecht2016spectral}. 
They are based on collocation methods which numerically solve a differential equation by selecting a finite number of collocation nodes over a time step, and approximate the solution on these nodes by a polynomial. 
A collocation method applied to a differential equation in time is in fact a subclass of fully-implicit Runge-Kutta methods, see~\cite{norsett1987solving}. 

SDC methods iteratively solve a collocation problem and can therefore be seen as iteratively converging to a Runge-Kutta scheme. 
The scheme begins from a low order solution on the collocation nodes; this solution is corrected with iteration ``sweeps" with the option to do any possible final update to calculate the solution at the next point in time. 
\cite{xia2007efficient} prove that each sweep increases the order of the scheme by one, up to the order of the collocation solution. 
Since the number of collocation nodes and correction sweeps can be supplied as runtime parameters, the scheme can be of arbitrary order without changing software infrastructure. 

A semi-implicit split version of SDC (SISDC) was proposed by~\cite{minion2003semi} for ODEs and later extended to incompressible flow, see~\cite{minion2004semi}. 
SISDC has been studied for a range of parabolic PDEs, where the stiff term is diffusive, for example for phase-field problems by~\cite{guo2019semi,liu2015stabilized}. 
Inspired by the study by~\cite{DurranBlossey2012} of IMEX multi-step methods for fast-wave-slow-wave problems (hyperbolic PDEs where the stiffness is due to a fast moving wave) a variant of SISDC called fast-wave slow-wave SDC (FWSW-SDC) was developed by~\cite{ruprecht2016spectral}.
They demonstrated that FWSW-SDC, in contrast to IMEX multi-step, is stable for arbitrarily large fast CFL numbers as long as the ``slow" time scale is resolved. 
In addition, they found that the computational cost of FWSW-SDC is comparable to IMEX Runge-Kutta schemes. 
In general, while FWSW-SDC costs more per timestep than an IMEX Runge-Kutta scheme, it is more stable and accurate and allows to take larger timesteps. 
Thus, if the number of required time steps is sufficiently reduced, the higher cost per step can be offset and SDC can be more efficient than Runge-Kutta methods in terms of work-precision.  
However, the examples for which FWSW-SDC was studied so far are relatively simple, linear and posed on a rectangular domain.
In this paper, we demonstrate that FWSW-SDC can also accurately solve more realistic benchmark problems, including advection on a sphere, nonlinear gravity wave in a channel, buoyancy-driven moist rising bubble and a dry baroclinic wave in a channel.

It was demonstrated by~\cite{jia2013spectral} that high-order fully implicit SDC methods perform well in ``climate-like" simulations. 
However, they use fully implicit SDC applied to the shallow water equations while we use the computationally more efficient FWSW-SDC. Both implicit SDC and FWSW-SDC require nonlinear solvers, however the nonlinear system to solve for FWSW-SDC is simpler, and therefore faster to solve.
We also extend these results by solving the compressible Euler equations, exploring the viability of FWSW-SDC methods to solve equations used by a modern dynamical core. 

FWSW-SDC could offer a promising alternative to the low order semi-implicit methods commonly used in atmospheric modelling.
It could be advantageous over the longer timescales seen in climate modelling. 
A second advantage of the iterative approach of SDC methods is that they allow for parallel\-isation-in-time to further increase their efficiency.
In the nomenclature introduced by~\cite{Gear1988}, parallel SDC offers parallelism-across-the-method, similar in spirit to parallel RKM with a (block-)diagonal Butcher table.
While this approach only offers limited concurrency, being able to use as many threads as there are collocation nodes, it can deliver speedup beyond saturation of spatial parallelism alone, see~\cite{FreeseEtAl2024}.
It also achieves better parallel efficiencies than larger-scale parallel-across-the-steps methods like Parareal by~\cite{LionsEtAl2001}, the parallel full approximation scheme in space and time (PFASST) by~\cite{emmett2012toward} or multigrid-reduction-in-time (MGRIT) by~\cite{FalgoutEtAl2014_MGRIT}, which can use a higher number of processing units in time.
Somewhat similar to parallel SDC are revisionist integral deferred correction (RIDC) by~\cite{ong2016algorithm}, which compute the corrections in a pipelined way meaning that each correction ``sweep" can be computed using a different processor.

Some meteorological centres, such as the Met Office, see ~\cite{walters2014met}, and the Deutscher Wetterdienst, see~\cite{zangl2015icon}, have unified modelling frameworks, where the same code base is used for both weather and climate modelling. 
A key benefit is that the maintenance of two models is not required, saving significant costs in staff and time. 
In addition,  if a model performs well across NWP and climate timescales, scientists can be more confident that developments to that model are better representing the atmosphere, see~\cite{walters2014met}. 
A disadvantage of this approach is that some modelling choices, such as order of numerical scheme, may be more suited to NWP than climate modelling. 
\cite{jia2013spectral} argue that high-order temporal accuracy may have significant benefits for climate modelling, owing to compounded temporal error over time. 
Conversely, the cost of using high order time discretisations for NWP may outweigh the benefits gained from the error reduction, so a mid-order scheme could be used. 
Having an arbitrary order scheme such as SDC allows for different order schemes for NWP and climate applications without violating the unified modelling framework.

For spatial discretization, we use a compatible finite element method described in Section~\ref{sec:model}\ref{subsec:spatial_disc}, comparable to that used in GungHo, the Met Office's next-generation dynamical core, see~\cite{melvin2019mixed, melvin2024mixed}. 
A key difference is that GungHo uses the lowest order compatible finite element spaces, whilst we used the next to lowest-order space. 
In addition, GungHo uses finite volume transport, whilst we use finite element transport.

The remainder of this paper is organised as follows. 
In Section~\ref{sec:model} we describe the governing equations, compatible finite element spatial discretisation, SDC time discretisation, and iterative solver. Section~\ref{sec:results} then analyses this formulation through a set of standard benchmark atmospheric modelling test cases for the compressible Euler equations.
\section{Model} \label{sec:model}
The formulation we present is implemented in Gusto, a dynamical core library that uses compatible finite elements (see \cite{bendall2020compatible} and \cite{hartney2024compatible} for other implementations using Gusto). It is built upon Firedrake, see \cite{FiredrakeUserManual}, a finite element library for solving partial differential equations. 
In particular, we make use of the extruded mesh capability by~\cite{Bercea2016}, the mesh management by~\cite{LangeMitchellKnepleyGorman2015}, automated symbolic manipulation of finite elements by~\cite{mcrae2016automated} and distributed parallel computing by~\cite{Dalcin2011}. 
Our code uses the `qmat' Python package from~\cite{lunet2024qmat} to generate the collocation nodes and time integration matrices, which are read in at the start of a run. 
The SDC algorithm is implemented in Gusto.
In this paper we use a method-of-lines approach where the spatial dimensions are first discretised, reducing the problem to an initial value problem which is then integrated in time with SDC.

\subsection{Governing Equations}
In this section we will describe the governing equations of our model. We solve the moist compressible Euler equations in the same form as \cite{bendall2020compatible} but without the rain moisture species. The prognostic variables are the wind velocity $\bm{u}$, the density of dry air $\rho$, the Exner pressure $\Pi$, the virtual dry potential temperature $\theta_{vd}$ and two moisture mixing ratios of water vapour $m_v$ and cloud water $m_c$. The virtual dry potential temperature $\theta_{vd}$ is defined as
\begin{equation}
     \theta_{vd} = T\left(\frac{p_R}{p}\right)^{R_d/c_{pd}}(1+R_v m_v/R_d) = \theta(1+R_v m_v/R_d) 
\end{equation}
where $T$ is the temperature, $p$ is the pressure, $c_{pd}$ is the specific heat capacity for dry air, $p_R$ is some reference pressure, $\theta$ is the dry potential temperature, $R_d$ is the specific gas constant for dry air, and $R_v$ is the specific gas constant of water vapour. The moist compressible Euler equation set is then
\begin{subequations} \label{eqn:compressible_euler}
\begin{align}
& \pfrac{\bm{u}}{t}+ \underbrace{(\bm{u} \cdot \nabla)\bm{u}}_\text{\clap{explicit~}} + \underbrace{\bm{f}\times\bm{u}}_\text{\clap{implicit~}} + \underbrace{\frac{c_{pd} \theta_{vd}}{1+m_t} \grad{\Pi}}_\text{\clap{implicit~}}+\underbrace{g\hat{\bm{k}}}_\text{\clap{implicit~}} = 0\label{eqn:compressible_euler_mom}, \\
& \pfrac{\rho}{t} + \underbrace{\rho \bm{\nabla} \cdot \bm{u}}_\text{\clap{implicit~}} + (\underbrace{\bm{u}_h \cdot \bm{\nabla}_h) \rho}_\text{\clap{explicit~}} + (\underbrace{\bm{u}_v \cdot \bm{\nabla}_v) \rho}_\text{\clap{implicit~}}= 0 \label{eqn:compressible_euler_rho}, \\
& \pfrac{\theta_{vd}}{t} + (\underbrace{\bm{u}_h \cdot \bm{\nabla}_h)\theta_{vd}}_\text{\clap{explicit~}} + (\underbrace{\bm{u}_v \cdot \bm{\nabla}_v)\theta_{vd}}_\text{\clap{implicit~}} + \underbrace{\mu \theta_{vd} \dv{\bm{u}}}_\text{\clap{implicit~}} = \underbrace{\mathcal{S}_{\theta_{vd}}}_\text{\clap{explicit~}}, \\ &
\pfrac{m_{v}}{t} + (\underbrace{\bm{u} \cdot \bm{\nabla})m_{v}}_\text{\clap{explicit~}}  = \underbrace{\mathcal{S}_{m_v}}_\text{\clap{explicit~}}, \\
& \pfrac{m_{r}}{t} +  (\underbrace{\bm{u} \cdot \bm{\nabla})m_{r}}_\text{\clap{explicit~}}  = \underbrace{\mathcal{S}_{m_r}}_\text{\clap{explicit~}}, 
\end{align}
\end{subequations}
where $g \hat{\bm{k}}$ is gravity and $\bm{f}=f\bm{\hat{k}}$ is the Coriolis parameter due to the planet's rotation multiplied by the unit vector in the vertical direction. $\bm{u}_v =  \hat{\bm{k}} (\hat{\bm{k}} \cdot \bm{u}$) is the vertical wind, and $\bm{u}_h = \bm{u} - \bm{u}_v$ is the horizontal wind. Similarly $\bm{\nabla}_v = \hat{\bm{k}} (\hat{\bm{k}} \cdot \bm{\nabla})$ is the vertical derivative and $\bm{\nabla}_h = \bm{\nabla} - \bm{\nabla}_v$ is the horizontal derivative.

The total moisture mixing ratio is defined as $m_t = m_v + m_c$, the factor $\mu$ is
\begin{equation}
    \mu = \left( \frac{R_m}{c_{vml}} - \frac{R_d c_{pml}}{c_{pd}c_{vml}} \right),
\end{equation}
and the source terms are
\begin{subequations} 
\begin{align}
S_{\theta_{vd}} & = \theta_{vd} \left[- \frac{c_{vd}L_V(T)}{c_{vml}c_{pd}T} - \frac{R_v}{c_{vml}}\left(1-\frac{R_dc_{vml}}{R_m c_{pd}} \right) -\frac{R_v}{R_m} \right]\Delta m_{c}^{cond},\\
S_{m_{v}} & = -\Delta m_{c}^{cond},\\
S_{m_{c}} & = \Delta m_{c}^{cond},\\
\end{align}
\end{subequations}
where $c_{vd}$, $c_{pml}$ and $c_{vml}$ are the specific heat capacities for dry air at a constant volume, moist air at a constant pressure and moist air at a constant volume. The moist air capacities $c_{pml}$ and $c_{vml}$ are functions of the mixing ratios.
In the absence of water, the factor $\mu$ is zero. $R_m$ is the specific gas constant and the latent heat from water vaporisation is given by $L_v(T)$, where $T$ is temperature. $\Delta m_{c}^{cond}$ represents the evaporation and condensation microphysics schemes described in \cite{bendall2020compatible}, and is a function of  $\rho$, $\theta_{vd}$, $m_v$ and $m_c$. 

The terms in \eqref{eqn:compressible_euler} are labeled explicit or implicit to indicate how they are integrated by the FWSW-SDC scheme defined in Section~\ref{sec:model}\ref{subsec:time_disc}. 
Transport terms have been treated explicitly, with the exceptions being the vertical advection terms for $\theta_{vd}$ and $\rho$, which control the fast gravity waves. The divergence terms $\rho \dv{\bm{u}}$ and $\mu \theta_{vd} \dv{\bm{u}}$ are associated with fast motions and are also treated implicitly.
The equation of state is
\begin{equation} \label{eqn:eos}
\Pi = \left(\frac{p}{p_R}\right)^{\kappa} = \left( \frac{\rho R_d \theta_{vd}}{p_R}\right)^{\kappa/(1-\kappa)}
\end{equation}
where $\kappa=R_d/c_{pd}$ and $\Pi$ is the Exner pressure.
\subsection{Spatial Discretisation} \label{subsec:spatial_disc}
In this section we describe the compatible finite element spatial discretisation used in our model, and the numerical advantages of using such a spatial discretisation. We will explain our choices of finite element function spaces, which will then allow us to describe the spatial discretisation of Equation \eqref{eqn:compressible_euler}. Our method-of-lines approach enables us to then describe the time discretisation separately in Section \ref{sec:model}\ref{subsec:time_disc}.

It has been demonstrated by~\cite{cotter2012mixed} that using families of finite elements that satisfy the discrete de Rahm complex can mimic the Arakawa C-grid staggering and hence satisfy many of the properties described in \cite{staniforth2012horizontal} without an orthogonal horizontal grid. In particular, this spatial discretisation provides good linear wave dispersion properties and avoids spurious pressure mode generation owing to factors such as the ratio between velocity and pressure degrees of freedom (DoFs).
These methods can be used with near uniform grids, avoiding the computational bottleneck at the poles seen with the orthogonal latitude-longitude mesh, see~\cite{cotter2012mixed}. The methods in \cite{cotter2012mixed} have been expanded to 3D Cartesian domains by \cite{natale2016compatible} and \cite{melvin2019mixed}. In particular, \cite{natale2016compatible} and \cite{melvin2018choice} presented function spaces choices to mimic the Charney-Phillips vertical grid staggering, allowing for benefits from improved vertical temperature gradient representation. This approach allows for arbitrarily high orders of accuracy whilst still maintaining these properties.

In 2D we use quadrilateral elements, and in 3D we use hexahedral elements with quadrilateral faces. In both cases the spaces are tensor products of horizontal and vertical spaces using the tensor product calculation of  \cite{mcrae2016automated}. 
We represent the prognostic variables by the following function spaces in each element
\begin{itemize}
    \item $\bm{u} \in \mathbb{V}_u$ is in the Raviart-Thomas function space of degree $1$, with continuous normal DoFs on element facets.
    \item $\rho \in  \mathbb{V}_{\rho}$ is the discontinuous linear scalar function space of degree $1$.
    \item $\theta_{vd}, m_v, m_c \in \mathbb{V}_{\theta}$ is quadratic continuous in the vertical and discontinuous linear in the horizontal.
\end{itemize}

\begin{table}[t] 
\centering
\begin{tabular}{c|c|c|c}  
\centering Space       & \centering $\mathbb{V}_u$ & \centering $\mathbb{V}_\theta$ & \centering $\mathbb{V}_{\rho}$ \tabularnewline
\hline
\centering Variables   & \centering $\bm{u}$ & \centering $\theta_{vd}$, $m_v$, $m_c$ & \centering $\rho$ \tabularnewline
\hline
2D    &  \centering \includegraphics[width=0.16\textwidth]{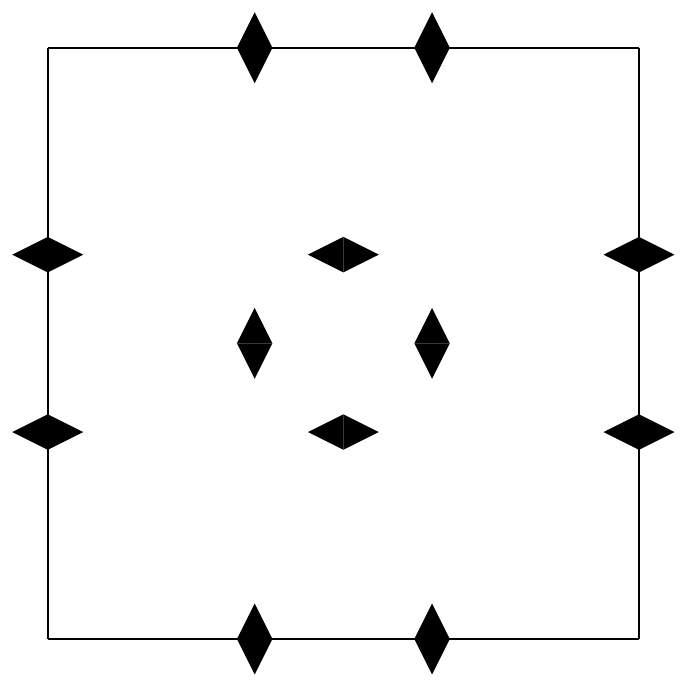} &
               \centering \includegraphics[width=0.16\textwidth]{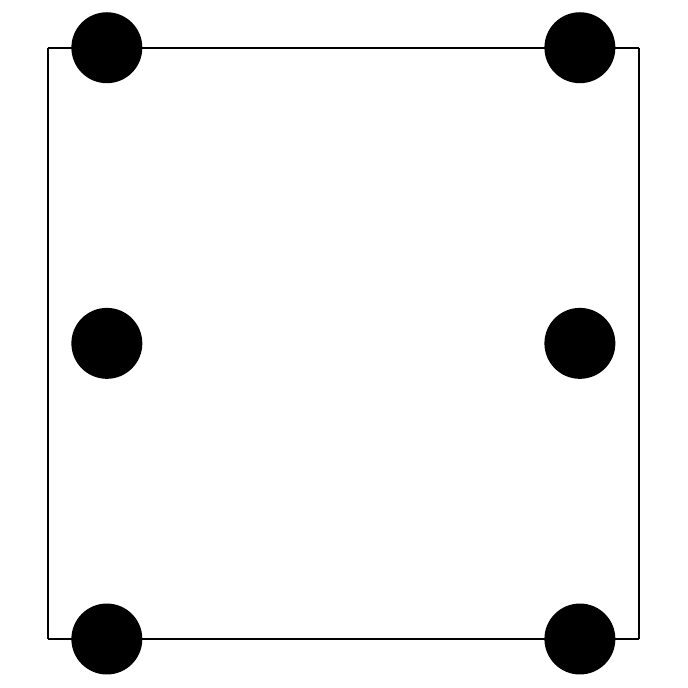} &
               \centering \includegraphics[width=0.16\textwidth]{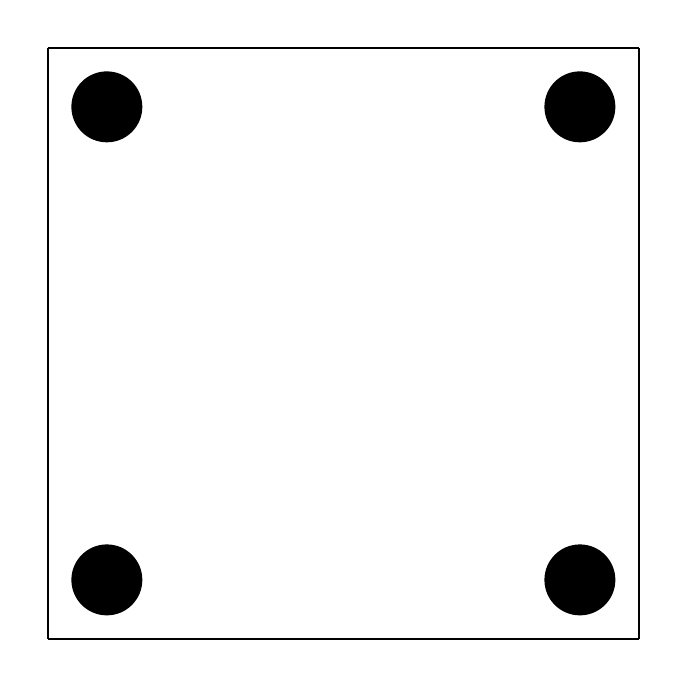} \tabularnewline
\hline
3D    &  \centering \includegraphics[width=0.16\textwidth]{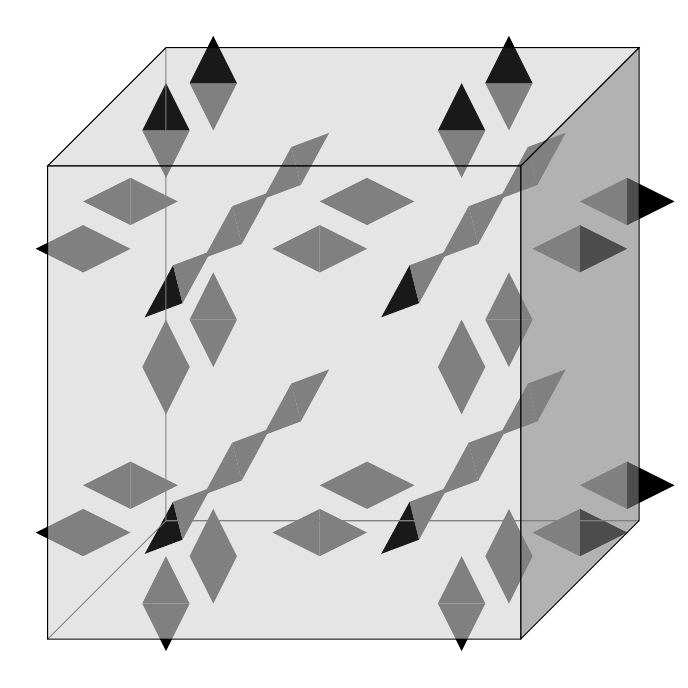} &
               \centering \includegraphics[width=0.16\textwidth]{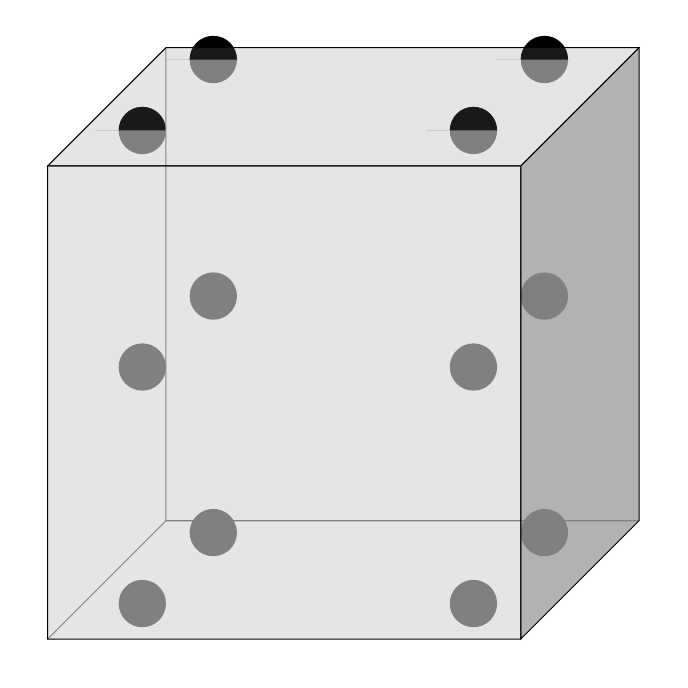} &
               \centering \includegraphics[width=0.16\textwidth]{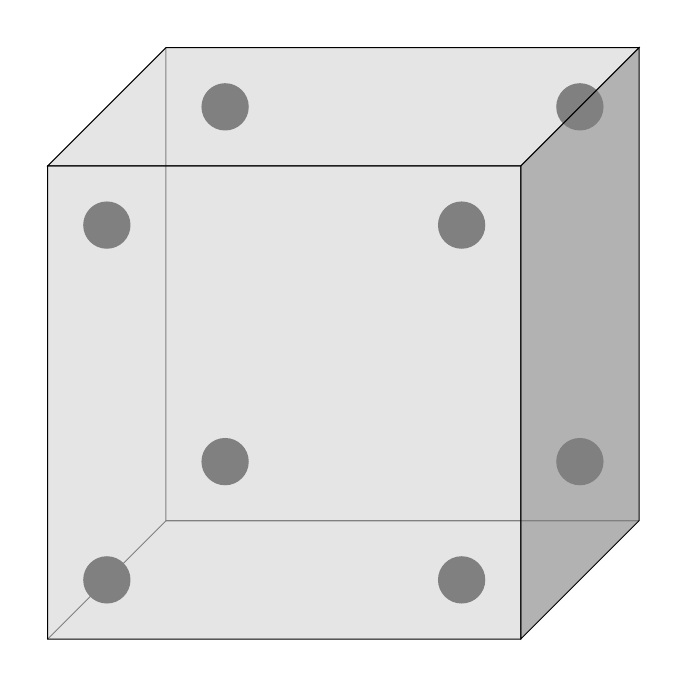} \tabularnewline
\end{tabular}
\caption{The finite elements function spaces used for wind, potential temperature, density and Exner pressure for the 2D and 3D problems explored. $\mathbb{V}_{u}$ and $\mathbb{V}_{\theta}$ have continuity between elements, hence sharing DoFs on element facets.}
\label{tab:elements}
\end{table}
These spaces are displayed in Table \ref{tab:elements}. The finite element formulation requires multiplication by the test functions $\bm{w} \in \mathbb{V}_u$, $\phi \in  \mathbb{V}_{\rho}$ and $\gamma \in \mathbb{V}_{\theta}$, followed by integration over the domain $\Omega$. The spatial discretisation of Equation \eqref{eqn:compressible_euler} is then

\begin{subequations} \label{eqn:spatial_discretisation}
    \begin{align}
        \begin{split}
                \int_{\Omega} \left( \bm{w} \cdot \partial_t \bm{u} + f\bm{w} \cdot \bm{\hat{k}} \times \bm{u} \right)dV \\ 
                -\int_{\Omega}  \bm{\nabla} \cdot(\bm{u} \otimes \bm{w} ) \cdot \bm{u} dV + \int_{\partial \Gamma} [[\bm{u} \otimes\bm{w}]]_{\bm{n}} \bm{u}^{*} dS \\
                -c_{pd}\int_{\Omega} \bm{\nabla} \cdot (\frac{\theta_{vd} \bm{w}}{1+m_t})\Pi dV + c_{pd}\int_{\Gamma} [[\frac{\theta_{vd} \bm{w}}{1+m_t}]]_{\bm{n}} \langle \Pi \rangle dS + \int_{\Omega} g(\hat{\bm{k}} \cdot \bm{w}) dV &= 0, \; \; \forall \bm{w} \in \mathbb{V}_u
        \end{split}
        \\
        \int_{\Omega} \left( \phi \partial_t \rho + \phi\rho\bm{\nabla}\cdot\bm{u} - \bm{\nabla}_h\cdot(\phi\bm{u}_h)\rho  - \bm{\nabla}_v\cdot(\phi\bm{u}_v)\rho \right) dV \\ +\int_{\Gamma} [[\phi \bm{u}_h]]_{\bm{n}} \rho^{*} + [[\phi \bm{u}_v]]_{\bm{n}} \langle \rho \rangle dS & = 0 , \; \; \forall \phi \in \mathbb{V}_{\rho}, \\
        \int_{\Omega} \left(\gamma \partial_t \theta_{vd}  + \mu\gamma\theta_{vd}\bm{\nabla}\cdot\bm{u} - \bm{\nabla}_h \cdot(\gamma^*\bm{u}_h)\theta_{vd} + \gamma^* \bm{u}_v \cdot \bm{\nabla}_v \theta_{vd} \right) dV 
        \\ + \int_{\Gamma_v} [[\gamma^* \bm{u}_h]]_{\bm{n}} \theta_{vd}^{*} dS & = \int_{\Omega} \gamma \mathcal{S}_{\theta_{vd}} dV, \; \; \forall \gamma \in \mathbb{V}_{\theta}, \\
        \int_{\Omega}\left( \gamma \partial_t m_{v} - \bm{\nabla} \cdot(\gamma^* \bm{u})m_{v} \right) dV +\int_{\Gamma_v} [[\gamma^* \bm{u}]]_{\bm{n}} m_{v}^{*} dS & = \int_{\Omega} \gamma \mathcal{S}_{m_{v}} dV , \; \; \forall \gamma \in \mathbb{V}_{\theta}, \\
        \int_{\Omega} \left(\gamma \partial_t m_{c} - \bm{\nabla} \cdot(\gamma^*\bm{u})m_{c}\right) dV +\int_{\Gamma_v} [[\gamma^* \bm{u}]]_{\bm{n}} m_{c}^{*} dS & = \int_{\Omega} \gamma \mathcal{S}_{m_{c}} dV, \; \; \forall \gamma \in \mathbb{V}_{\theta}.
    \end{align}
\end{subequations}
Here, $dV$ corresponds to integrating over an element and $dS$ to element facets. $\Omega$ is the domain, $\Gamma$ is the set of all interior facets, $\partial \Gamma$ is the external domain boundary and $\Gamma_v$ is the set of all vertical facets between columns. 
We have integrated some terms in the equation by parts where the continuity of the field means that the differential operator is not defined, which results in the surface terms. 
The notation $\partial_t$ represents the partial derivative  $\frac{\partial}{\partial t}$ with respect to time $t$. $[[a]]$ is a jump $a^{+}-a^{-}$ for a scalar $a$, where the superscripts $+$ and $-$ refer to variable values either side of an element facet; 
$[[\bm{x}]]_{\bm{n}}$ is $(\bm{x}^{+}-\bm{x}^{-})\cdot \bm{n}$ for a vector $\bm{x}$; 
$\bm{x} \otimes \bm{y}$ is the outer product for vectors $\bm{x}$ and $\bm{y}$ and $[[\bm{x}\otimes \bm{y}]]_n = (\bm{n}^+ \cdot \bm{x}^+)\bm{y}^+ + (\bm{n}^- \cdot \bm{x}^-)\bm{y}^-$. $\langle a \rangle$ is an average  $(a^{+}+a^{-})/2$  across an element facet; 
$\bm{x}^*$ and $a^*$ correspond to upwind values of $\bm{x}$ and $a$ on a element facet, choosing $\bm{x}^{+}$ and $a^{+}$ or $\bm{x}^{-}$ and $a^{-}$ dependent on whether $\bm{u}\cdot \bm{n} \geq 0$. For transforms between physical space and reference element, $\mathbb{V}_u$ functions use the contravariant Piola transform, whilst $\mathbb{V}_{\rho}$ and $\mathbb{V}_{\theta}$ functions have no explicit mapping. Further details are available in~\cite{mcrae2016automated}. 

We use the Streamline Upwind Petrov–Galerkin (SUPG) stabilization method by~\cite{brooks1982streamline, tezduyar1988petrov, tezduyar1989finite} in the vertical for fields in $\mathbb{V}_{\theta}$. This reduces spurious oscillations associated with continuous finite element transport schemes. This alters the test function to $\gamma^* \rightarrow \gamma + \tau \bm{k} \cdot \bm{u} \partial_{z}\gamma$ which provides an upwind bias. Here, $\tau$ is a stabilisation parameter chosen to be $\tau = 1/\sqrt{15}$. We only alter the test function in the transport terms, not the time derivative terms, as this provided the most stable results.

\subsection{Time Discretisation} \label{subsec:time_disc}
In this section we write a single timestep of standard SDC as a preconditioned unmodified Richardson iteration, following the approach in \cite{huang2006accelerating, ruprecht2016spectral,weiser2018theoretically,speck2019algorithm}, applied to a general ordinary differential equation in time
\begin{equation} \label{eqn:ODE}
\bm{x}_t(t) = \bm{f}(\bm{x}(t)) \; \; \forall t.
\end{equation}
with initial conditions $\bm{x}(t_n)=\bm{x}_n.$
We then explain how to modify this to FWSW-SDC.
This will enable us to introduce the terminology required to describe our approach without the unnecessary complexity of the terms arising from the spatial discretisation, which are absorbed into the term $\bm{f}(\bm{x}(t))$. $\bm{x}$ is a vector of prognostic variables $\bm{x}=[\bm{u}, \rho, \theta_{vd},m_v]^T$ containing all the DoFs in the domain.

Equation \eqref{eqn:ODE} can be written in Picard integral form
\begin{equation} \label{eqn:int}
    \bm{x}(t) = \bm{x}_n + \int_{t_n}^{t}\bm{f}(\bm{x}(s)) ds,
\end{equation}
which we can then solve using collocation, integrating Lagrange polynomials across the $M$ collocation nodes $\tau_1,...,\tau_M$ where $t_n \leq \tau_1 \le ... \le \tau_M \leq t_{n+1}$. 
A standard choice for collocation nodes are nodes such as Gauss-Legendre, Gauss-Radau or Gauss-Lobatto. The SDC scheme can be of order up to the underlying quadrature rule, hence $2M$, $2M-1$ and $2M-2$ for Gauss-Legendre, Gauss-Radau and Gauss-Lobatto respectively.
We then approximate Equation \eqref{eqn:int} as
\begin{equation} \label{eqn:quad}
    \bm{x}_{m} = \bm{x}_n + \Delta t \sum_{j=1}^M q_{m,j}\bm{f}(\bm{x}_{j}),
\end{equation}
since
\begin{equation} \label{eqn:int2}
    \int_{t_n}^{\tau_{m}}\bm{f}(\bm{x}(s)) ds \approx \Delta t\sum_{j=1}^M q_{m,j}\bm{f}(\bm{x}_{j}),
\end{equation}
where $\bm{x}_m \approx \bm{x}(\tau_m)$.
The $q_{m,j}$ are quadrature weights on the interval $[t_n, \tau_m]$, which we calculate by integrating Lagrange polynomials.
\begin{equation} \label{eqn:lagrange}
 q_{m,j} = \int_{t_n}^{\tau_{m}} l_j(s) ds
\end{equation}
\begin{equation} \label{eqn:lagrange1}
    l_j(t) = \frac{1}{c_j} \prod_{k=1, k \neq j}^{M} (t - t_k), \; \;  c_j = \prod_{k=1, k \neq j}^{M} (t_j - t_k).
\end{equation}
We write this as a system of equations in matrix form, resulting in the collocation problem
\begin{equation}
    (\bm{I} - \Delta t\bm{Q}\bm{f})(\bm{X})=\bm{X}_n,
\end{equation}
where $\bm{X} = [\bm{x}_1, \bm{x}_2, \cdots, \bm{x}_M]^T \in \mathbb{R}^{M}$ is a vector containing our solution on each collocation node, $\bm{X}_n = [\bm{x}_n, \bm{x}_n, \cdots, \bm{x}_n]^T \in \mathbb{R}^{M}$ is a vector containing $M$ copies of the initial conditions, $\bm{Q} = [q_{i,j}]_{i,j} \in \mathbb{R}^{M\times M}$ is an integration matrix, where each entry $q_{i,j}$ is a quadrature weight calculated by integrating Lagrange polynomials. Finally, $\bm{f}(\bm{X}) = [\bm{f}(\bm{x}_1), \bm{f}(\bm{x}_2), ..., \bm{f}(\bm{x}_M)]^T  \in \mathbb{R}^{M}$ contains the non-time derivative terms on each collocation node. As discussed before, this is similar to a fully-connected implicit method such as implicit Runge-Kutta. It is potentially a large problem, with $M$ copies of the model state. In order to solve more efficiently we use a Richardson Iteration, which results in
\begin{equation}
    \bm{X}^{k+1}=\bm{X}^k + (\bm{X}_n - (\bm{\bm{I}}- \Delta t \bm{Q}\bm{f})(\bm{X}^{k}))
\end{equation}
where each correction sweep is $k=1,...,K$. However, this would be relatively slow to converge. Preconditioning using a lower triangular quadrature matrix $\bm{Q}_{\Delta}$, that approximates $\bm{Q}$, improves convergence. 
Using the preconditioner $(\bm{I} -\Delta t \bm{Q}_{\Delta}\bm{f})$ results in
\begin{equation}
 \bm{X}^{k+1}=\bm{X}^k + (\bm{I} -\Delta t \bm{Q}_{\Delta}\bm{f})^{-1}(\bm{X}_n - (\bm{\bm{I}}- \Delta t \bm{Q}\bm{f})(\bm{X}^{k})).
\end{equation}
Multiplying by $(\bm{I} -\Delta t \bm{Q}_{\Delta}\bm{f})$ results in
\begin{equation}
 (\bm{I} -\Delta t \bm{Q}_{\Delta}\bm{f})\bm{X}^{k+1}=(\bm{I} -\Delta t \bm{Q}_{\Delta}\bm{f})\bm{X}^k + (\bm{X}_n - (\bm{\bm{I}}- \Delta t \bm{Q}\bm{f})(\bm{X}^{k}))
\end{equation}
By expanding and canceling out terms we get
\begin{equation}
 (\bm{I} -\Delta t \bm{Q}_{\Delta}\bm{f})(\bm{X}^{k+1})=\bm{X}_n + \Delta t (\bm{Q}-\bm{Q}_{\Delta})\bm{f}(\bm{X}^{k}).
\end{equation}
This is the \textit{zero-to-node} formulation of SDC. equivalent to the Gauss-Seidel method. For each collocation node, an SDC update is
\begin{equation}
\bm{x}_{m}^{k+1} - \Delta t \hat{q}_{m,m}\bm{f}(\bm{x}^{k+1}_m) = \bm{x}_n + \Delta t\sum_{j=1}^{m-1} \hat{q}_{m,j} \bm{f}(\bm{x}^{k+1}_j) + \Delta t\sum_{j=1}^{M} (q_{m,j}-\hat{q}_{m,j})\bm{f}(\bm{x}^{k}_j)
\end{equation}
where each $\hat{q}_{i,j}$ is an entry in the matrix $\bm{Q}_{\Delta}=[\hat{q}_{i,j}]_{i,j}$. If $\bm{Q}_{\Delta}$ is strictly lower triangular, the $q_{m,m}$ term on the left hand side vanishes we would have an explicit method. The final collocation update to time level $n+1$ is then
\begin{equation}  \label{eqn:final_update}
\bm{x}_{n+1} = \bm{x}_n + \Delta t\sum_{j=1}^{M} w_{m,j}\bm{f}(\bm{x}^{K}_j),
\end{equation}
where
\begin{equation} \label{eqn:lagrange_fin}
 w_{j} = \int_{t_n}^{t_{n+1}} l_j(s) ds.
\end{equation}
Equation \eqref{eqn:final_update} improves the order of accuracy by $1$ if the scheme has not reached the maximum order of the underlying collocation problem, see~\cite{ruprecht2016spectral}. 

For the fast-wave slow-wave (FWSW) method of \cite{ruprecht2016spectral}, we split $\bm{f}=\bm{F}+\bm{S}$, where $\bm{F}$ represents terms responsible for fast motions, and $\bm{S}$ represents terms responsible for slow motions. We also choose both an implicit and an explicit $\bm{Q}_{\Delta}$. In this paper we use the explicit Euler and \verb|MIN-SR-NS| by \cite{CaklovicEtAl2025} as $\bm{Q}_{\Delta}$ matrices for the slow motions and the implicit Euler, the so called LU-trick by \cite{Weiser2015}, and the \verb|MIN-SR-FLEX| by~\cite{CaklovicEtAl2025}, $\bm{Q}_{\Delta}$ matrices for the fast motions. They are defined below.

Implicit Euler and explicit Euler preconditioners are given by 
\begin{equation}
\bm{Q}^{\verb|IE|}_{\Delta} = \frac{1}{\Delta t}\left(\begin{array}{cccc}
\Delta \tau_1 & 0 & \cdots & 0\\
\Delta  \tau_1 & \Delta \tau_2 & \cdots & 0\\
\vdots & \vdots & \ddots & 0\\
\Delta \tau_1 & \Delta \tau_2 & \cdots & \Delta \tau_M
\end{array} \right)\qquad
\bm{Q}^{\verb|EE|}_{\Delta} =\frac{1}{\Delta t} \left(\begin{array}{cccc}
0  & \cdots & \cdots & 0\\
\Delta \tau_2 & 0& \cdots & 0\\
\vdots & \ddots & \ddots & 0\\
\Delta \tau_2 & \cdots & \Delta \tau_M & 0
\end{array} \right),\qquad
\end{equation}
respectively. For the LU trick, we take the LU-decomposition of $\bm{Q}^T$, and choose $\bm{Q}^{\verb|LU|}_{\Delta} = \bm{U}^T$. This is a common choice in the SDC literature, owing to its improved convergence over the implicit Euler $\bm{Q}^{\verb|IE|}_{\Delta}$. The \verb|MIN-SR-NS| and the \verb|MIN-SR-FLEX| precondtioners are given by $\bm{Q}^{\verb|MIN-SR-NS|}_{\Delta} = \diag(\tau_i/M)$ and $\bm{Q}^{\verb|MIN-SR-FLEX|}_{\Delta} = \diag(\tau_i/k)$ with $k \in \{1, \dots, M\}$. They minimise the spectral radius of the so called iteration matrix in the non-stiff and the stiff limit, respectively, and, hence, improve convergence. Both preconditioners allow SDC to run in parallel due to their diagonal structure, however, in this work we run the method serially. Note that $\bm{Q}^{\verb|MIN-SR-FLEX|}_{\Delta}$ changes with each iteration. For details we refer to \cite{CaklovicEtAl2025}.
\\
\\ 
The \textit{zero-to-node} FWSW-SDC formulation we use in this paper is
\begin{equation}
\begin{split}
(1 - \Delta tq^{imp}_{m,m}\bm{F})(\bm{x}_{m}^{k+1}) = \bm{x}_n + \Delta t\sum_{j=1}^{m-1}(q^{imp}_{m,j} \bm{F}(\bm{x}^{k+1}_j) + q^{exp}_{m,j} \bm{S}(\bm{x}^{k+1}_j)) \\ + \Delta t\sum_{j=1}^{M} ((q_{m,j} - q^{imp}_{m,j})\bm{F}(\bm{x}^{k}_j) +(q_{m,j} - q^{exp}_{m,j}) \bm{S}(\bm{x}^{k}_j)),
\end{split}
\end{equation}
where $q^{imp}_{i,j}$ is an entry in the implicit matrix $\bm{Q}^{Imp}_{\Delta}=[q^{imp}_{i,j}]_{i,j}$ and $q^{exp}_{i,j}$ is an entry in the lower-triangular matrix $\bm{Q}^{Exp}_{\Delta}=[q^{exp}_{i,j}]_{i,j}$. A nonlinear solve is required here for the implicit terms, but FWSW SDC treats some terms explicitly, speeding up convergence of the nonlinear solve.

Our initial solution is $\bm{x}_n=[\bm{u}, \rho, \theta_{vd}, m_v, m_c]_n^T$, and our spatially discretised non-time derivative terms are then
\begin{subequations}
\begin{align}
\bm{F}(\bm{x}) &= \begin{bmatrix}
            \int_{\Omega} f\bm{w} \cdot \bm{\hat{k}} \times \bm{u}dV- c_{pd}\int_{\Omega} \bm{\nabla} \cdot \left( (\theta_{vd} \bm{w})/(1+m_t) \right)\Pi dV \\ + c_{pd}\int_{\Gamma} [[(\theta_{vd} \bm{w})/(1+m_t)]]_{\bm{n}} \langle \Pi \rangle dS + \int_{\Omega} g(\hat{\bm{k}} \cdot \bm{w}) dV, \\
            \int_{\Omega} \phi\rho\bm{\nabla}\cdot\bm{u} -  \bm{\nabla}_v\cdot(\phi\bm{u}_v)\rho  dV +\int_{\Gamma} [[\phi \bm{u}_v]]_{\bm{n}} \langle \rho \rangle dS, \\
            \int_{\Omega} \mu\gamma\theta_{vd}\bm{\nabla}\cdot\bm{u} + \gamma^* \bm{u}_v \cdot \bm{\nabla}_v \theta_{vd} dV,            \\
            0,
            \\
            0.
           \end{bmatrix} \\
\bm{S(\bm{x})} &= \begin{bmatrix}
            - \int_{\Omega}  \bm{\nabla} \cdot(\bm{u} \otimes \bm{w} ) \cdot \bm{u} dV + \int_{\partial \Gamma} [[\bm{u} \otimes\bm{w}]]_{\bm{n}} \bm{u}^{*} dS, \\
            -  \int_{\Omega}  \bm{\nabla}_h\cdot(\phi\bm{u}_h)\rho dV +  \int_{\Gamma} [[\phi \bm{u}_h]]_{\bm{n}} \rho_{d}^{*} dS,\\
           - \int_{\Omega}  \bm{\nabla}_h \cdot(\gamma^*\bm{u}_h)\theta_{vd} dV +\int_{\Gamma_v} [[\gamma^* \bm{u}_h]]_{\bm{n}} \theta_{vd}^{*} dS - \int_{\Omega} \gamma \mathcal{S}_{\theta_{vd}} dV, \\
           - \int_{\Omega}  \bm{\nabla} \cdot(\gamma^*\bm{u})m_{v} dV +\int_{\Gamma_v} [[\gamma^* \bm{u}]]_{\bm{n}} m_{v}^{*} dS - \int_{\Omega} \gamma \mathcal{S}_{m_{v}} dV,
           \\
           - \int_{\Omega}  \bm{\nabla} \cdot(\gamma^*\bm{u})m_{c} dV +\int_{\Gamma_v} [[\gamma^* \bm{u}]]_{\bm{n}} m_{c}^{*} dS - \int_{\Omega} \gamma \mathcal{S}_{m_{c}} dV.
           \end{bmatrix}
\end{align}
\end{subequations}

The SDC Algorithm \ref{alg:SDC} describes the solution procedure, where we first calculate an initial state on the collocation nodes, either by copying across the initial state, or using a low order timestepping method. Next we compute successive corrections across the nodes, each correction adding up to an order or accuracy, up to the order of the underlying collocation method. Finally, a final collocation step computes the solution at the next time level, this also adds up to an order of accuracy, see~\cite{ruprecht2016spectral}.
\begin{algorithm}
\caption{Spectral deferred correction algorithm pseudocode. At the start, either a low order guess is computed on the collocation nodes using an IMEX Euler scheme, or the initial state at time level $n$ is just copied across to the nodes. Successive correction updates are then computed, with a nonlinear implicit solve. Finally, a collocation update is computed, or if the collocation nodes are such that the final node is at the next time level, the final node value can be copied.}
\begin{algorithmic}
\For{n = 1,\;N}  \Comment{Timestep loop}
\State $\bm{x}^{0}_{0} \gets \bm{x}_n$
    \For{m = 1,\;M}
        \Comment{Compute first guess on collocation nodes}
        \If{low order guess}
            \State $\bm{x}^{0}_m \gets \mathrm{Solve}(\bm{x} -  \Delta t\bm{F}(\bm{x}) = \bm{x}^{0}_{m-1} + \Delta t\bm{S}(\bm{x}^{0}_{m-1})) \; \mathrm{for} \; \bm{x}$
        \ElsIf{copy initial state}
            \State $\bm{x}^{0}_m \gets \bm{x}_n$
        \EndIf
    \EndFor
    \For{k = 1,\;K} \Comment{Loop through K correction iterations}
        \For{m = 1,\;M} \Comment{Compute nodal correction}
        \State \(
         \begin{aligned}[t]
            \bm{R}_m \gets 
            &\Delta t\sum_{j=1}^{m-1}(q^{imp}_{m,j} \bm{F}(\bm{x}^{k+1}_j) + q^{exp}_{m,j} \bm{S}(\bm{x}^{k+1}_j)) + \notag \\ 
            &\Delta t\sum_{j=1}^{M} ((q_{m,j} - q^{imp}_{m,j})\bm{F}(\bm{x}^{k}_j) +(q_{m,j} - q^{exp}_{m,j}) \bm{S}(\bm{x}^{k}_j)) \notag
           \end{aligned} \)
        \State $\bm{x}^{k+1}_m \gets \mathrm{Solve}(\bm{x} -  \Delta tq^{imp}_{m,m}\bm{F}(\bm{x}) = \bm{x}_n + \bm{R}_m) \; \mathrm{for} \; \bm{x}$

        \EndFor
    \EndFor
    \If{final collocation update} \Comment{Compute final update}
            \State $\bm{x}_{n+1} \gets \bm{x}_{n} + \Delta t\sum_{j=1}^{M} w_j(\bm{F}(\bm{x}^{K}_{j})+\bm{S}(\bm{x}^{K}_{j}))$
    \ElsIf{copy final node value}
            \State $\bm{x}_{n+1} \gets \bm{x}^{K}_{M}$
    \EndIf
\EndFor
\end{algorithmic}
\label{alg:SDC}
\end{algorithm}

The choice of node type, whether a low order scheme is used and whether a final update is computed will be stated in the description of the test cases in Section \ref{sec:results}. If Gauss-Legendre nodes are used, then the final collocation update step is required, but if Gauss-Radau or Gauss-Lobatto nodes are used, then the final collocation update is not required.
\subsection{Iterative Solver}
In order to efficiently solve the compressible Euler equations an SDC method, we require an effective solution method for the nonlinear system arising in steps 2 and 3 of Algorithm \ref{alg:SDC}. We achieve this through the following strategy, which makes use of the PETSc linear algebra library by~\cite{petsc-efficient} \cite{petsc-user-ref} that is used by Firedrake and therefore Gusto. We take a similar approach to \cite{cotter2023compatible}, which makes use of PETSc's Newton linesearch method, facilitated by Firedrake's  \texttt{NonlinearVariationalSolver} object. As discussed in \cite{cotter2023compatible}, the linear solve then uses a GMRES method, preconditioned with an Additive Schwartz Method (ASM). The ASM preconditioner does a direct vertical solve in star patch subdomains using a LU factorisation.

We found that our model still converges with the expected accuracy and had no stability issues using slacker solver tolerances than in \cite{cotter2023compatible}. GungHo (see~\cite{melvin2019mixed, melvin2024mixed}) uses comparable tolerances to good effect.
The nonlinear solver's absolute and relative tolerances are $tol^{nl}_a = tol^{nl}_r = 10^{-4} $, and the linear solver absolute and relative tolerances are $tol^{l}_a = tol^{l}_r = 10^{-4}$. In addition, we use the method by~\cite{eisenstat1996choosing} to dynamically adjust the linear solver tolerances. The result is the nonlinear solver typically converging in $1-3$ iterations, and the linear solver converging in $4-10$ iterations. We need $(M)\times(K)$ implicit solves where $M$ is the number of collocation nodes and $K$ the number of SDC iterations.

\section{Results} \label{sec:results}
In this section, we demonstrate how FWSW-SDC can be used effectively in atmospheric modelling by exploring a variety of NWP test cases. First we test the temporal convergence of the SDC schemes using an advection test case of \cite{williamson1992standard}. We use the non-hydrostatic gravity wave test case of \cite{skamarock1994efficiency} to explore the generation of gravity waves, and test the spatial and temporal convergence of the model. The moist rising bubble test case of \cite{bryan2002benchmark} tests the model with buoyancy driven motion and the evolution of moisture species. Finally, the dry baroclinic wave in the channel test case of \cite{ullrich2015analytical} measures the models ability to predict synoptic scale weather patterns. For the dry test cases in this paper, we do not use the moisture mixing ratio fields.

We use the notation $SDC(M,K)$ to describe SDC schemes, where $M$ is the number of collocation nodes and $K$ is the number of correction iterations. In general we use Gauss-Legendre nodes, copy the initial conditions and do a final collocation update, therefore the order of the scheme should be $K+1$. 

\subsection{Advection on the sphere} \label{subsec:adv}
The first test case is an advection test case of \cite{williamson1992standard}. We use it as a self-convergence test for SDC schemes, in order to eliminate the spatial error which often dominates the solution. We advect a tracer $D$ satisfying 
\begin{equation}\label{eqn:adv} 
    \pfrac{D}{t} + \bm{u} \cdot \grad{D} = 0.
\end{equation}
The initial condition of the tracer is
\begin{equation}
D = \left\lbrace
\begin{matrix}
 \frac{1}{2}D_{max} (1+ \cos{\left(\frac{3 \pi r}{R}\right)}) & \mathrm{if} \: r \: \le \: R, \\
0 & \mathrm{otherwise},
\end{matrix} \right.
\end{equation}
where $D_{max}=1000$m, $r$ is the great circle distance between the longitude and latitude $(\lambda, \phi)$ and  $(\lambda_c, \phi_c)=(3 \pi/2, 0)$. $R=a/3$, where $a=6.37122 \times 10^6$m is the radius of the Earth. The prescribed wind velocity is taken as the perpendicular gradient of the streamfunction:
\begin{equation}
    \psi = -a u_{max}(\sin{\phi}\cos{\alpha}-\cos{\lambda}\cos{\phi}\sin{\alpha}) 
\end{equation}
where here we set $\alpha=0$, $u_{max}=2 \pi a/(12 \mathrm{days})$.

We solve Equation \eqref{eqn:adv} on the sphere with a C$32$ mesh, where C$n$ refers to a cubed-sphere mesh with $n \times n$ horizontal cells on each of the $6$ cubed-sphere panels. We use $SDC(2,3)$, $SDC(3,5)$, $SDC(4,7)$ with Gauss-Legendre nodes. The scheme is explicit and uses $Q^{\verb|EE|}_{\Delta}$, copying initial conditions and computing the final collocation update. We use $\Delta t = [2400, 1800, 1200, 900]$s, and the ``true" solution is calculated using a SSPRK3 time scheme with $\Delta t_{true}=0.5$s. We use a conjugate gradient linear solver with a block Jacobi preconditioner. Each block is preconditioned using ILU. The linear solver is required to invert the mass matrices.

\begin{figure}[t]
\centering
  \centering
  \includegraphics[width=.7\linewidth]{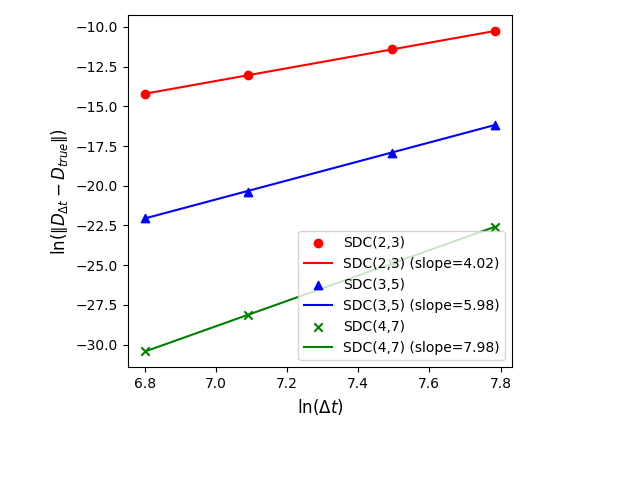}
  \caption{Temporal convergence of SDC for an advection test on the sphere. We use $\Delta t = [2400, 1800, 1200, 900]s$ on a C$32$ cubed sphere mesh. We used a self convergence test to eliminate the spatial error with a SSPRK3 scheme and $\Delta t_{true} = 0.5$s.  $SDC(2,3)$ and $SDC(3,5)$ achieve the expected order of convergence of $4$, $6$ and $8$ respectively.  }
\label{fig:time_conv}
\end{figure}

As can be seen from Figure \ref{fig:time_conv}, when calculating the normalised $L_2$ error the expected convergence is achieved for $SDC(2,3)$, $SDC(3,5)$ and $SDC(4,7)$. $SDC(2,3)$ should be $4$th order, $SDC(3,5)$ $6$th order and $SDC(4,7)$ $8$th order. 
\subsection{Non-hydrostatic gravity waves} \label{subsec:nh_gw}
The non-hydrostatic gravity wave test case of \cite{skamarock1994efficiency} models generation of gravity waves that then propagate through a 2D vertical slice domain of size  $(L,H) = (300, 10)$km. An initial potential temperature perturbation of the form
\begin{equation}
    \theta^{'} = \Delta \theta_0 \frac{\sin(\pi z/H)}{1+(x-x_c)^2/a^2}
\end{equation}
is added to a background state. Here $\Delta \theta_0 = 10^{-2}$K, $a=5$km and $x_c=0$m. The background potential temperature is specified with a constant buoyancy frequency $\mathcal{N}=0.01\;\mathrm{s}^{-1}$. $\rho$ is then calculated to be in hydrostatic balance. The background wind is $U=20\mathrm{ms}^{-1}$.  This perturbation then spreads out forming inertia-gravity waves.

As discussed in \cite{melvin2010inherently}, the original test case is for the Boussinesq equations, whilst \cite{hundertmark2007regularization} updated it for the compressible Euler equations. 
However, the results are largely similar. Here we solve the compressible Euler equations.

For the solution we use $SDC(2,3)$ with Gauss-Legendre nodes, LU implicit $Q_{\Delta}$ and forward Euler explicit $Q_{\Delta}$. A final collocation update is completed, but the initial guess on the nodes is just a copy of the start of timestep state. This scheme is $4$th order in time.

\begin{figure}[t]
  \centering
  \includegraphics[width=\linewidth]{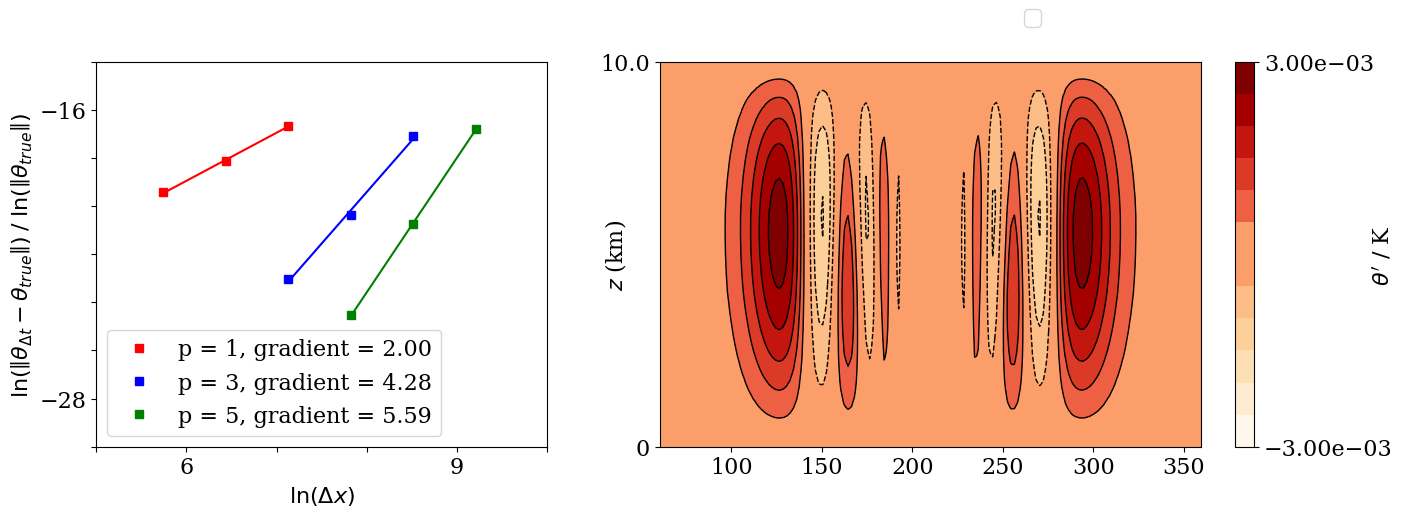}
  \caption{(Right) Final $\theta^{'}$solution for the gravity-wave test case at $t=3000$s for $(\Delta x, \Delta z) =(2000, 1000)$m and $\Delta t = 6$s. Contours are spaced every $5 \times 10^{-4}$K. The solution is the same as in the literature \cite{melvin2010inherently}, \cite{skamarock1994efficiency}. (Left) Convergence of $\theta$ $L_2$ error for the gravity-wave test case at $t=3000$s for a constant Courant–Friedrichs–Lewy (CFL) number of $0.06$ for the $p = 1$ case. The $p = 1$ $SDC(2,3)$ scheme has second order convergence, $p = 3$ $SDC(3,5)$ has better than fourth order convergence and $p = 5$ $SDC(4,7)$ has better than fifth order convergence.}
\label{fig:gw_sol}
\end{figure}

The final solution at $t=3000$s seen in Figure \ref{fig:gw_sol} (Right) for $(\Delta x, \Delta z) =(2000, 1000)$m and $\Delta t = 6$s is in line with literature. In addition, we ran a self-convergence test for various spatial orders $p = 1$, $p = 3$ and $p = 5$, where $p$ refers to the order of the discontinuous space in the compatible finite element function spaces. For the $p = 1$, $p = 3$ and $p = 5$ spatial schemes we use $SDC(2,3)$, $SDC(3,5)$ and $SDC(4,7)$ time discretisations with Gauss-Radau collocation nodes, LU implicit $Q_{\Delta}$ and Forward Euler explicit $Q_{\Delta}$. No final collocation update is done, and the initial guess on nodes is once again a copy of the start of timestep state. These schemes should be $3$rd, $5$th and $7$th order in time. 

For all schemes we have $\Delta z = 1000$m, and the ``true" solution, calculated from high resolution simulations, have $\Delta t_{true} = 0.15$s. Table \ref{tab:resolution} displays the resolutions used for the convergence test for each of the schemes. 

\begin{table}[t] 
\centering
\begin{tabular}{|c|c|c|c|c|}
\hline
$p$ & $SDC$ scheme & $\Delta x$ (m) & $\Delta t$ (s) &  $\Delta x_{true}$ (m) \\
\hline
1 & $SDC(2,3)$ & $[1250, 625, 312.5]$ & $[3.75, 1.875, 0.9375]$ & $50$ \\
\hline
3 & $SDC(3,5)$ & $[5000, 2500, 1250]$ & $[3.75, 1.875, 0.9375]$ & $200$\\
\hline
5 & $SDC(4,7)$ & $[10000, 5000, 2500]$ & $[1.875, 0.9375, 0.46875]$ & $800$\\
\hline
\end{tabular}

\caption{Table summarizing the grid resolutions, time steps sizes and resolution of ``true" solution for different order schemes.}
\label{tab:resolution}
\end{table}

The $L_2$ error against ``true" solution was calculated at the native resolutions using the point data interpolation of \cite{nixon2024consistent}. Figure \ref{fig:gw_sol} (Left) shows that the scheme converges at $2$nd order for $p = 1$ $SDC(2,3)$, better than $4$th order for $p = 3$ $SDC(3,5)$ and better than $5$th order for $p = 5$ $SDC(4,7)$. This result demonstrates that using an SDC timestepping scheme with a finite element spatial discretisation can provide methods with variable and arbitrarily high order of accuracy, defined at runtime.

\subsection{Moist rising bubble} \label{subsec:bubble}
The moist rising bubble test case of \cite{bryan2002benchmark} simulates a buoyancy-driven thermal bubble through a saturated atmosphere in a vertical slice domain. As the bubble rises, it overturns symmetrically on each side. The domain is of dimensions $(L,H) = (10, 10)$km, where the horizontal boundaries are periodic and the vertical boundaries are zero-flux. We use the microphysics scheme as described in \cite{bendall2020compatible}. Here, as in \cite{bryan2002benchmark} and \cite{bendall2020compatible}, there is no Coriolis force.

The initial conditions follow the procedure of \cite{bendall2020compatible} to approximately satisfy the following background state:
\begin{itemize}
    \item Constant total moisture mixing ratio $m_t = 0.02\mathrm{kgkg}^{-1}$,
    \item Constant wet equivalent potential temperature $\theta_e = 320$K,
    \item Water vapour mixing ratio at saturation,
    \item Surface pressure of $p=10^5$Pa,
    \item Hydrostatically balanced,
    \item Zero initial velocity.
\end{itemize}
The following perturbation is then added to the virtual dry potential temperature $\theta_{vd}$
\begin{equation} \label{eqn:conservation_config}
\theta^{'}_{vd} = \left\lbrace
\begin{matrix}
 \Delta \Theta \cos^2{\left(\frac{\pi r}{2r_c}\right)} & \mathrm{if} \: r \: \le \: c_r, \\
0 & \mathrm{otherwise},
\end{matrix} \right .
\end{equation}
where $r = \sqrt{(x-x_c)^2+(z-z_c)^2}$, $\Delta \Theta=2K$, $r_c=2$km, and $(x_c,z_c)=(L/2,2)$km. As in \cite{bendall2020compatible} the $\theta_{vd}$ field is $\theta_{vd} = \bar{\theta}_{vd}(1+\theta^{'}_{vd}/300K)$.
Following~\cite{bendall2020compatible}, $\rho$ is found such that the pressure field is unchanged by the perturbation, and $m_v$ is calculated such that it is at saturation.

As in Section \ref{sec:results}\ref{subsec:nh_gw}, we use the $4$th order $SDC(2,3)$ with Gauss-Legendre nodes, copying the initial conditions to get the initial guess on the nodes and completing the final collocation update to get the end of timestep solution. We use two formulations, one using $Q^{\verb|LU|}_{\Delta}$ and $Q^{\verb|EE|}_{\Delta}$ and the other with $Q^{\verb|MIN-SR-FLEX|}_{\Delta}$ and $Q^{\verb|MIN-SR-NS|}_{\Delta}$.

\begin{figure}[t]
    \centering

    \includegraphics[width=\textwidth]{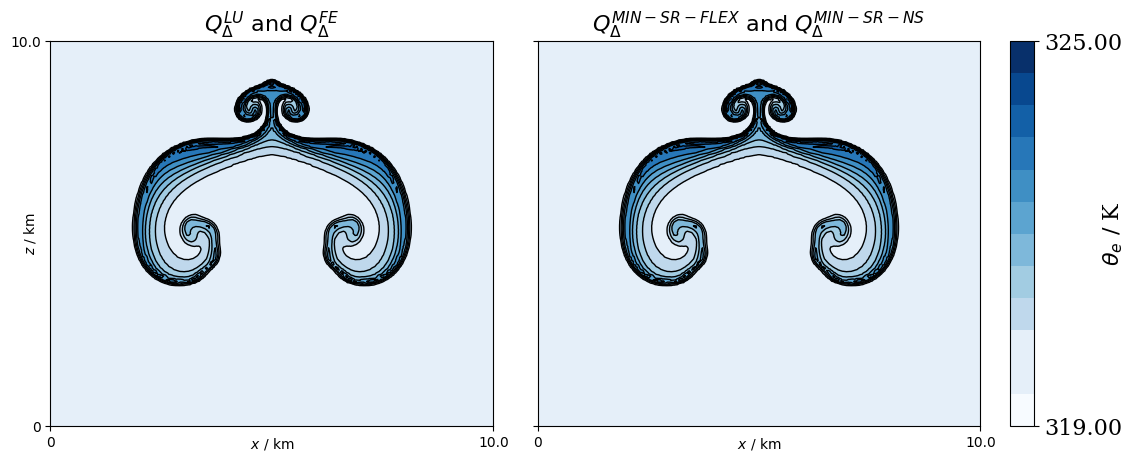} 
    \caption{$\theta_e$ field at time $t = 1000$ s for $Q^{\texttt{LU}}_{\Delta}$ and $Q^{\texttt{ExpEuler}}_{\Delta}$ (Left), and $Q^{\texttt{MIN-SR-FLEX}}_{\Delta}$ and $Q^{\texttt{MIN-SR-NS}}_{\Delta}$ (Right). As in \cite{bendall2020compatible}, $\Delta x = \Delta z = 100$m and $\Delta t = 1$ s. The $320$K contour has been omitted, and the contours are spaced every $0.5$K.}

    \label{fig:moist_bf}
\end{figure}
From Figure \ref{fig:moist_bf} we see that for both formulations the results replicate the next to lowest-order results seen in \cite{bendall2020compatible} at $t=1000$s with the same resolution, where the second plume owing to a Rayleigh-Taylor instability also appears. The diagonal preconditioners $Q^{\verb|MIN-SR-FLEX|}_{\Delta}$ and $Q^{\verb|MIN-SR-NS|}_{\Delta}$ produce very similar results to the commonly chosen $Q^{\verb|LU|}_{\Delta}$ and $Q^{\verb|EE|}_{\Delta}$, but have the possibility of time parallelisation. 

\subsection{Baroclinic wave} \label{subsec:baro}
The dry baroclinic wave test case of \cite{ullrich2015analytical} models the generation of synoptic-scale waves induced by baroclinic instability in a $3$D Cartesian channel with Coriolis effects. The channel has the dimensions $(L_x, L_y, H) = (40000, 6000, 30)$km for the $(x,y,z)$ directions, with periodic boundaries in the $x$ direction and zero-flux boundary conditions in the $(y,z)$ directions.

\begin{figure}[t]
    \centering
    \includegraphics[width=\textwidth]{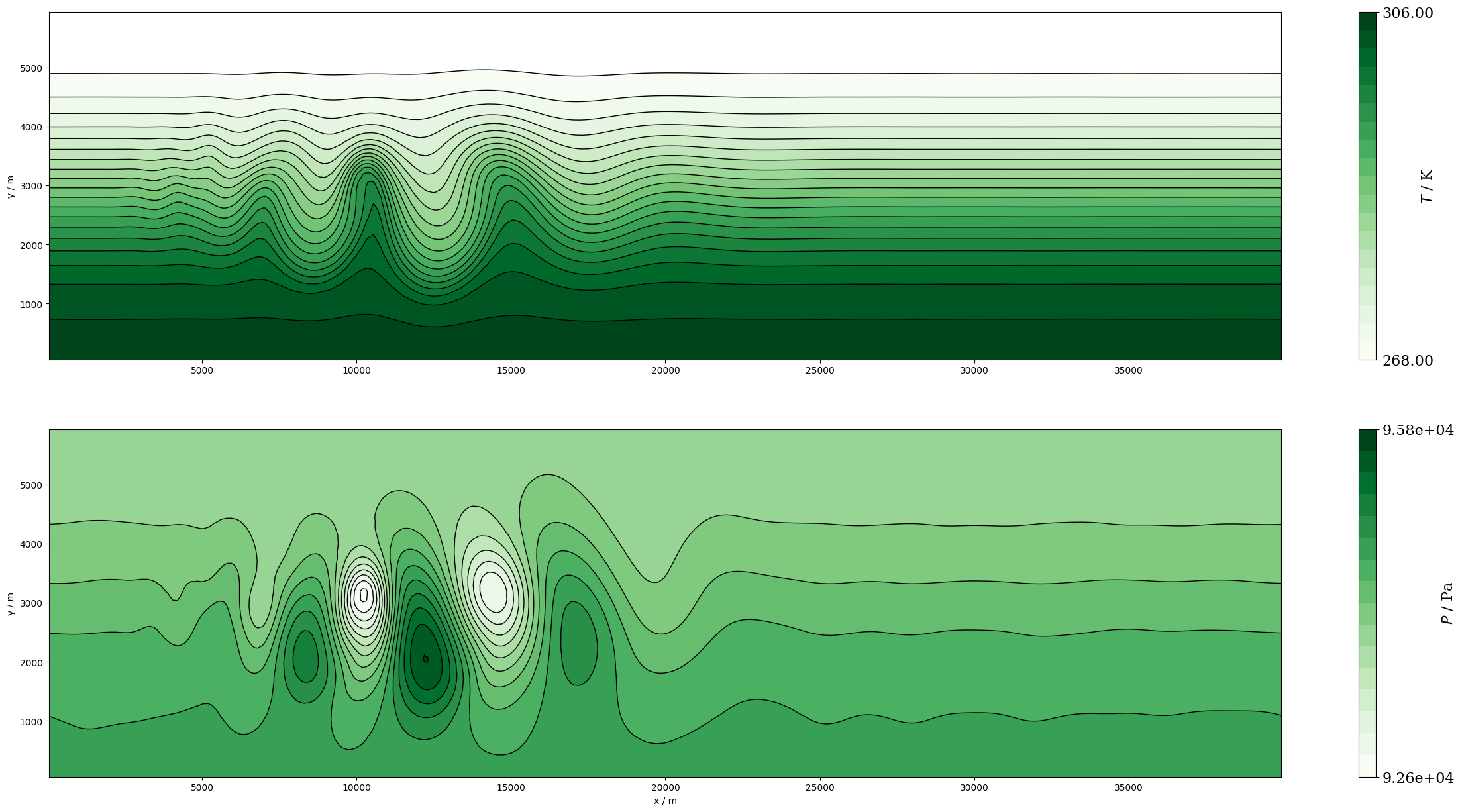} 
    \caption{Temperature with contours from $268$K to $306$K with spacings of $2$K (top) and pressure with contours from $92600$Pa to $95800$Pa with spacings of $200$Pa (bottom) fields at $t=12$days $z=500$m, the resolutions were $\Delta t = 1800$s, $\Delta x = \Delta y = 250$km and $\Delta z = 1.5$km. }
    \label{fig:baroclinic}
\end{figure}

The initial background flow field is a mid-latitude zonal jet in thermal wind balance. The initial conditions for the zonal wind $\bm{u}$, the temperature $T$ and the geopotential $\Phi$ are therefore defined as

\begin{subequations} \label{eqn:baroclinic}
\begin{align}
&    \bm{u}(x,y,z) = -\bm{u}_0 \sin^2 {\left (\frac{\pi y}{L_y} \right )}\ln{\eta} \exp{ \left (-\left(\frac{\ln{\eta}}{b}\right)^2\right)} \\
& T(x,y,\eta) = T_0\eta^{R_d\Gamma/g} + \frac{f_0 u_0}{2 R_d} \left[ y - \frac{L_y}{2} - \frac{L_y}{2\pi}\sin{\left(\frac{2\pi y}{L_y}\right)}  \right] \left[ \frac{2}{b^2}(\ln{\eta})^2 \right] \exp{\left(-\left(\frac{\ln{\eta}}{b}\right)^2\right)}  \\
& \Phi(x,y,\eta) = \frac{T_0 g}{\Gamma}\left(1- \eta^{ R_d\Gamma/g}\right) +  \frac{f_0 u_0}{2 }\left[ y - \frac{L_y}{2} - \frac{L_y}{2\pi}\sin{\left(\frac{2\pi y}{L_y}\right)}  \right]\ln{\eta} \exp{\left(-\left(\frac{\ln{\eta}}{b}\right)^2\right)}  
\end{align}
\end{subequations}
where $\bm{u}_0=35$ms$^{-1}$, $T_0=288$K,$f_0=2.0 \time 10^{-6}$Km$^{-1}$ $b=2$ and $g = 9.80616$ms$^{-2}$. The vertical pressure coordinate $\eta$ is used in \cite{ullrich2015analytical}. As in \cite{bendall2020compatible} the initial conditions are converted to $(u, \rho, \theta)$, finding $\eta$ through Newton iteration. We use $\Phi=gz$ with $\eta \in \mathbb{V}_{\theta}$, and solve for $\eta$ with
\begin{equation}
    \eta^{k+1} = \eta^{k} - \frac{\Phi^{k} - g}{T^{k}_v - R_d/\eta^{k}}
\end{equation}
$\eta$ is then used to calculate $\bm{u}$, $T$ and $\Phi$. Unlike in \cite{bendall2020compatible}, we have no moisture in this test case. We therefore calculate $\theta_{vd}$ as
\begin{equation}
  \theta_{vd} = T\left(\frac{1}{\eta}\right)^{R_d/c_{pd}}
\end{equation}
The dry density $\rho$ is calculated to be in hydrostatic balance. Finally, a perturbation is added to the background velocity to induce baroclinic instability. For $(x_c, y_c)=(2000, 2500)$km, $L_p = 6 \times 10^6$ we use
\begin{subequations}
\begin{align}
    & r(x,y) = \sqrt{(x - x_c)^2 - (y - y_x)^2} \\
    & u^{'}(x,y) = \exp{\left(-\left(\frac{r}{L_p}\right)^2\right)} 
\end{align}
\end{subequations}
We use the $4$th order $SDC(2,3)$ with Gauss-Legendre nodes. Figure \ref{fig:baroclinic} shows the results from a $12$ day simulation with a timestep of $\Delta t = 1800$s, $\Delta x = \Delta y = 250$km and $\Delta z = 1.5$km.  
The baroclinic wave develops as in~\cite{ullrich2015analytical}, with similar areas of low and high pressure. The later developing regions of high and low pressure are slightly weaker, however this is likely due to the resolution difference between the simulations.

\section{Conclusion and Outlook}
This work has expanded on previous SDC literature, such as~\cite{jia2013spectral} who applied implicit SDC to the shallow water equations and~\cite{ruprecht2016spectral} who applied FWSW-SDC to linear problems, exploring convergence and stability. Here we have applied the FWSW-SDC variant to classic atmospheric dynamics test cases, in order to explore the viability of using SDC as a time discretisation in NWP and climate models. A key advantage of SDC over schemes such as Runge-Kutta is the flexibility gained from the arbitrary order.
As mentioned in the introduction, this could be particularly useful for the unified modelling approach for NWP and climate models that is often adopted at meteorological organisations.

In Section~\ref{sec:results}\ref{subsec:adv} we demonstrated that the SDC scheme matches expected convergence up to order $8$ for an advection problem on the sphere. We then explored the convergence order of schemes with differing spatial and temporal orders, with the arbitrary order of our spatial and temporal schemes. Here we achieved $2$nd, better than $4$th and better than $5$th order convergence in Section~\ref{sec:results}\ref{subsec:nh_gw}, showing that the combination of SDC in time and finite elements in space can achieve variable and arbitrary high order of accuracy. 

A moist rising bubble test case was used in Section \ref{sec:results}\ref{subsec:bubble}, we compared $LU$-trick and explicit Euler SDC preconditioners against diagonal (parallelisable) preconditioners from~\cite{CaklovicEtAl2025}, showing almost identical results, comparable to works by~\cite{bryan2002benchmark, bendall2020compatible}. Time-parallel SDC will be explored in future work, but the diagonal preconditioners are robust and stable and the small-scale parallel speedup they can offer could give SDC an advantage over other schemes.

Finally, in Section \ref{sec:results}\ref{subsec:baro} we simulated a baroclinic wave in a Cartesian domain, a classic NWP test case for simulating synoptic-scale weather systems. SDC was stable for time steps similar to those used in \cite{ullrich2015analytical} and produced comparable results.

SDC offers a promising alternative to IMEX Runge-Kutta and other semi-implicit temporal discretisations used in NWP and climate modelling. The method's stability, arbitrary accuracy and substantial literature of time-parallel adaptations leave potential for it to be competitive. In this paper, we have explored its viability when running classic NWP and climate test case. Future work will look at time-parallel versions of SDC.

\section*{Acknowledgments}
We acknowledge funding from the Advanced Parallel in Time Algorithms for Partial Differential Equations (APinTA PDEs) project (grant: SPF EX20-8 Exposing Parallelism: Parallel-in-Time (DN517492)). JF and DR thankfully acknowledge funding from the German Federal Ministry of Education and Research (BMBF) under grant 16ME0679K. Supported by the European Union - NextGenerationEU.
For the purpose of open access, the author has applied a Creative Commons Attribution (CC BY) licence to any Author Accepted Manuscript version arising from this submission.

\section*{Data Statement}
The results in this paper were generated using Gusto, a dynamical core toolkit based on Firedrake. The source code is available on GitHub (\href{https://github.com/firedrakeproject/gusto}{\texttt{https://github.com/firedrakeproject/gusto}}). The scripts to run the tests, generate the results, and plot them are available in a public GitHub repository (\href{https://github.com/atb1995/brown_sdc_2025_scripts}{\texttt{https://github.com/atb1995/brown\_sdc\_2025\_scripts}}). The plotting scripts are dependent on \texttt{tomplot} (\href{https://github.com/tommbendall/tomplot}{\texttt{https://github.com/tommbendall/tomplot}}).
\bibliography{bibliography}
\end{document}